\documentclass{elsarticle}
\usepackage{graphicx} 
\usepackage{amsmath}
\usepackage{amssymb}
\usepackage{mathrsfs}
\usepackage{amsthm,tikz}
\usepackage{comment}
\usepackage{ulem}
\usepackage[colorlinks=true]{hyperref}
\usepackage{booktabs,multirow}
\usepackage{subcaption}
\usepackage{todonotes}
\newcommand{\abs}[1]{\left\lvert#1\right\rvert}
\newcommand{\udef}{\mathrel{\mathop:}=}

\newcommand{\R}{\mathbb{R}}

\renewcommand{\min}{{\mathrm{min}}}
\renewcommand{\max}{{\mathrm{max}}}

\newcommand{\correct}[1]{{\color{blue}{#1}}}
\newcommand{\first}[1]{{\color{black}{#1}}}
\newcommand{\second}[1]{{\color{black}{#1}}}

\theoremstyle{remark}
\newtheorem{remark}{Remark}
\newtheorem{definition}{Definition}
\newtheorem{prop}{Proposition}
\date{}

\journal{MATCOM}
\begin{document}

\begin{frontmatter}
\title{State Dependent Riccati for dynamic boundary control to optimize irrigation in Richards' Equation framework}
\author[first]{Alessandro Alla}
\affiliation[first]{organization={Sapienza, Università di Roma, Italy (corresponding author)},email={alessandro.alla@uniroma1.it}}
 \author[second]{Marco Berardi} 
\affiliation[second]{organization={Consiglio Nazionale delle Ricerche, Istituto di Ricerca sulle Acque, Italy}, email={marco.berardi@cnr.it}}
\author[third]{Luca Saluzzi}
\affiliation[third]{organization={Centro di Ricerche Matematiche "Ennio De Giorgi", Scuola Normale Superiore, Italy},email={luca.saluzzi@sns.it}}
\begin{abstract}
We present an approach for the optimization of irrigation in a Richards' equation framework. We introduce a proper cost functional, aimed at minimizing the amount of water provided by irrigation, at the same time maximizing the root water uptake, which is modeled by a sink term in the continuity equation. The control is acting on the boundary of the dynamics and
due to the nature of the mathematical problem we use a \first{State Dependent} Riccati approach which provides suboptimal control in feedback form, applied to the system of ODEs resulting from the Richards' equation semidiscretization in space.
The problem is tested with existing hydraulic parameters, also considering proper root water uptake functions. The numerical simulations also consider the presence of noise in the model to further validate the use of a feedback control approach.
\end{abstract}


\end{frontmatter}
\section{Introduction}\label{sec:intro}

It has been assessed that almost the 70\% of freshwater withdrawals is used in the agricultural sector, of which 86\% for irrigation. In this context, optimizing water needed for irrigation, without reducing the yield of the crops, becomes crucial for a sustainable management of agricultural practices: in this framework, our paper aims at providing efficient modeling and numerical tools for optimizing irrigation, resorting to control techniques. As a matter of facts, more comprehensive irrigation models rely on water flow in partially saturated soils. 

Recently, some control theory approaches have been used for optimizing irrigation processes. In \cite{MLNLD18}, for example, a linear parameter varying (LPV) \correct{model is developed for controller design}, aimed at maintaining soil water content in the root zone within a certain target. A very simplified model is proposed in \cite{CC19}, where a constant diffusivity is imposed to Richards' equation, and a sliding mode control is presented, whereas an optimal control approach is used in \cite{LFPDG16} for a \correct{basic} water balance law. And yet, a \correct{streamlined} tool is introduced in \cite{Berardi_et_al_TiPM_2022}, based on the computation of steady solutions of Richards’ equation. An interesting approach for applying control techniques in a Richards’ equation framework is described in \cite{WCISA19}, yet
with very different applications, i.e. maximizing
the amount of absorbed liquid by redistributing the materials,
when designing the material properties of a diaper. 

In \cite{BDG23}, it has been proposed a model for solving an optimal control problem under the quasi-unsaturated assumptions, which provides a suitable hydrological setting 
that prevents to reach water moisture saturation in the soil: the authors have derived the appropriate optimality conditions for the boundary control of a class of nonlinear Richards’ equations, and implemented these results in the development and computation of numerical solutions by a classical Projected Gradient Descent algorithm: here we will follow some concepts of that paper, with respect to the cost functional, albeit in a completely different control technique and therefore also  numerical framework. 
\medskip

In this work, we present an infinite horizon dynamic boundary control problem applied to the optimization of the irrigation. Loosely speaking, the idea of a dynamic boundary control problem is to include an additional dynamic on the boundary, which is modified according to the control under consideration. Indeed, dynamic boundary control involves the control of the derivative of the state equation at the boundary. This technique is commonly used to influence the evolution of partial differential equations governing the state of the system. 
This kind of problem has not been studied before to the best of our knowledge in the context of Richards' equation. A dynamic boundary control problem has been studied for finite horizon problems in e.g. \cite{KN04,KMMR22} under other circumstances as heat equation and models with drift terms. 
The objective of our control procedure is to maximize the root water uptake, at the same time minimizing the amount of water provided for irrigation, which is described by properly assigning the boundary conditions. To this purpose, a feedback control is implemented by a tailored numerical method.
Our goal is to obtain a reliable numerical method that allows us to compute the control in feedback form. This control is usually more suitable for application driven problem since it is able to react to perturbations or measurement errors which may occur very often when working with real life problems. In this work, we use the so-called {\em discretize-than-optimize} approach. Unlike the approach adopted in \cite{BDG23}, we initially discretize the problem using finite difference methods, resulting in a system of $d$ ordinary differential equations. Note that the dynamic boundary control is directly included in the discretized system. An optimal control in feedback form can be obtained through the knowledge of the value function, but its computation suffers from the curse of the dimensionality. Indeed, it is possible to compute the value function by means of the Hamilton-Jacobi-Bellman (HJB, see \cite{BCD97}) equation. However, this equation is not easy to solve either analytically, due to its nonlinear nature, or numerically, since the dimension of the problem can be \second{arbitrarily} high. The term {\em curse of dimensionality} has been coined due to the fact that increasing the dimension of the problem its complexity increases exponentially. Recently, the mitigation of the curse of dimensionality has been studied by different points of view. Among all, we mention: max-plus algebra \cite{Akian_Gaubert_Lakhoua_2008,maxplusdarbon}, sparse grids and polynomials \cite{GK16,AKK21,KVW23}, tree-structure algorithms \cite{AKS19,alla2019high,falcone2023approximation}, deep neural networks \cite{Darbon_Langlois_Meng_2020,Kunisch_Walter_2021,Zhou_2021,Onken2021,ruthotto2020machine,ABK21}, low-rank tensor decompositions \cite{richter2021solving,DKK21,DKS23} and kernel interpolation techniques \cite{alla2021hjb,ehring2023hermite}.

A viable alternative to compute an approximation of the value function is the use of the State Dependent Riccati Equation \correct{(SDRE)}, see e.g. \cite{BLT07,BTB00}. In this context, it is assumed that the value function can be locally approximated by a quadratic form. Furthermore, the state equation needs to have a peculiar form which has to be linear for a given state. Thus, the control can be computed through the solution of several algebraic Riccati equations, i.e. $d-$dimensional matrix equations. This approach has been studied theoretically in \cite{BH18}, and computationally variants to the method have been introduced in \cite{AKS23}. A model reduction approach to the topic can be found in \cite{HW23,HKW24,AP24}.\\
To summarize, the novelties presented in this work are (i) introduction of a dynamic boundary control problem for optimizing irrigation in a Richards' equation framework, (ii) the use of \correct{SDRE} to obtain the control in feedback form and (iii) simulations affected by the presence of noise.\\
The outline of the paper is the following. In Section \ref{sec:ocp}, we introduce the Richards' equation for modeling unsaturated flow with root water uptake macroscopic effects, and we present our dynamic boundary control problem. Later in Section \ref{sec:num}, we explain our discretization approach and the SDRE applied to our control problem. Numerical tests are discussed in Section \ref{sec:test}, referring to standard root water uptake functions and well-known hydraulic models and parameters. Finally, we highlight the novelties of this work, drawing some future research directions in this framework.

\section{Richards' equation and the control problem}\label{sec:ocp}
The so-called Richards' equation (RE) is the standard tool for modeling water flow in unsaturated soils: such flow is driven by capillarity and gravity forces. It has been derived by the combination of Darcy's law for porous media and mass balance equation, and it is the core of many disciplines in porous media, ranging from soil science, agronomy, environmental engineering, hydrogeology, numerical analysis and applied mathematics. For instance, the significant review \cite{Paniconi_Putti_WRR_2015} provides a perspective of the RE in the context of catchment hydrology, while the encyclopedic work in \cite{Vereecken_et_al_VZJ_2016} takes into account the importance of such model in soil science scientific community. From a mathematical point of view, since exact solutions of RE have been found only for specific settings (e.g. \cite{Tracy_WRR_2006}), efforts have been concentrated mainly in the numerical solvers. Many numerical issues can be considered when solving such equation: for instance, spatial discretization can be accomplished by finite volumes and mixed finite elements methods, as in \cite{Eymard_Comput_Geosci_1999,Schneid_Knabner_Radu_2004,Kees_Farthing_Dawson_CMAME_2008,Radu_Pop_Knabner_2008}, techniques for handling nonlinearities have been widely treated in  \cite{List_Radu_2016,Casulli_Zanolli_SIAM_2010,Berardi_Difonzo_Notarnicola_Vurro_APNUM_2019,Berardi_Difonzo_Lopez_CAMWA_2020}. By the way, the interested reader is referred to \cite{Farthing_Ogden_SSSAJ_2017} for a comprehensive review on numerical methods for RE.


Soil moisture dynamics is primarily governed by gravitational forces. \first{Therefore, for many agro-hydrological applications, studying vertical infiltration alone can provide valuable insights. However, significant heterogeneities and uncertainties in model parameters can sometimes render the use of Richards' equation across extensive two-dimensional spatial domains impractical. In line with previous agronomic applications, including those at the scale of irrigation districts (e.g., \cite{Hassan_et_al_Irr_Sci_2024, Coppola_et_al_Eco_2024}), we adopt a simplified approach by focusing on the one-dimensional, vertical formulation of Richards' equation, specifically in its pressure-based form,  having the following form:}

\begin{equation}
\label{eq:Richards}
C(h)\frac{\partial h}{\partial t}=\frac{\partial}{\partial z}\Big( K(h)\Big(\frac{\partial h}{\partial z}-1\Big)\Big)-S(h), \quad (z,t)\in [0,L]\times [0,T].
\end{equation}
Here, $h$ is the soil metric potential ($\rm cm$),  $C(h)$ is the soil water capacity ($\rm cm^{-1}$), $K(h)$ is the unsaturated hydraulic conductivity ($\rm cm \;d^{-1}$), $L$ is depth ($\rm cm$) of the domain, $T$ is the period of observation of the process. The sink term $S(h)$ represents the macroscopic water uptake by roots: for the sake of simplicity, it is assumed it depends only on the state variable $h$. \\
There exist different empirical  forms for $K(h)$ and $C(h),$ where the latter is defined as 
\begin{equation}\label{eq:c}
C(h) := \frac{\partial\theta}{\partial h},
\end{equation}
and where the {\em water retention curve} $\theta(h)$ has a suitable form (see, for example, \cite{Haverkamp,VanGenutchen}). \first{More details about the physical meaning of these functions are thoroughly described in the inspiring books  \cite{Nimmo_Enc,Bear_Cheng}: loosely speaking, while the most basic measure of water in porous media is the volumetric water content $\theta$, defined as the volume of water per bulk volume of the medium.
The matric pressure $h$, from now on simply pressure, can be thought of as the pressure difference across the water-air interface. Lastly, the water retention curve arises from the nonlinear fitting of laboratory measurements of $\theta$ and $h$ according to given empirical functions.}

Since RE is the core of most sophisticated irrigation models, here we want to provide some tools enhancing the  trade-off between the physical process and the anthropogenic action to sustainably manage water resources in irrigation framework: in practice the goal is to maximizing the root water uptake while at the same time minimizing water volumes provided by irrigation\correct{. This} task is accomplished by properly controlling the upper boundary condition according to a dynamic boundary control problem.

Let us now introduce the dynamic boundary control problem (\cite{KN04,KMMR22}) we are going to study in the present work:
\begin{align}\label{eq:state}
\begin{aligned}
C(h)\frac{\partial h}{\partial t}&=\frac{\partial}{\partial z}\Big( K(h)\Big(\frac{\partial h}{\partial z}-1\Big)\Big)-S(h), \quad (z,t)\in [0,Z]\times [0,T],\\
h(0,0) &= h_T , \\
\frac{\partial h}{\partial t}(0,t) &= u(t),\quad t\in[0,T],\\
h(Z,t) & =  h_B(t), \quad t\in[0,T],\\
h(z,0) & = h_0(z), \quad z\in[0,Z],
\end{aligned}
\end{align}

where $u(t):[0,T]\rightarrow \mathbb{R}$ is the control, $h_T \in \mathbb{R}$ is the initial condition for the upper boundary of the domain and $h_B:[0,T]\rightarrow \mathbb{R}$ is a given function acting on the lower boundary of the domain. We remark that usually, for boundary control systems there is no additional dynamic on the boundary as we are introducing in \eqref{eq:state}. This problem has been studied in \cite{KN04} for heat equation and recently in \cite{KMMR22} for equations involving drift terms. 
\first{As the control directly influences the boundary dynamics, the spatial semidiscretization of the system naturally assumes a form well-suited for applying the SDRE method. Specifically, the resulting system takes the form $y'(t)=A(y(t))y(t)+B u(t)$, where $y(t)$ represents the spatial semidiscretization of the solution $h(z,t)$.}

We are interested in minimizing the following tracking problem for the uptake function:
\begin{equation}\label{eq:cost}
    J(u) = \int_0^T \int_0^Z \left(|S(h(z,t)) - S_{max}|^2 + \lambda |u(t)|^2\, \right)dt \,dz
\end{equation}
where $S_{max}= \max_{h \in \mathbb{R}} S(h)$ and $\lambda>0$ a given parameter. Finally, our optimal \second{control problem} reads
\begin{equation}\label{ocp}
    \min_u J(u) \mbox{ such that } h \mbox{ solves } \eqref{eq:state}
\end{equation}

Existence and uniqueness results for the boundary control problem, can be found in \cite{Berardi_Difonzo_Guglielmi_2023}. \first{The extension to dynamic boundary control in this context is open and goes beyond the scopes of this work. We remark that there is a difference between the setting proposed in \cite{Berardi_Difonzo_Guglielmi_2023} which does not allow to extend their results straightforward.}

\section{Numerical approximation of the control problem}\label{sec:num}

In this work we use the {\em discretize-than-optimize} approach. Specifically, we first semidiscretize the dynamical system \eqref{eq:state} and then we apply a suboptimal feedback control strategy based on \correct{SDREs} (see e.g. \cite{BTB00,BLT07}).

First of all, we discretize the state equation \eqref{eq:state} using Finite Differences (FD, \cite{L07}). We consider a spatial mesh constituted by $d+1$ nodes equally spaced with stepsize $\Delta z$ such that $z_{i}=i \Delta z\ge 0$. Thus, we are looking for the vector function $y:[0,T]\to \R^{d}$ such that
\begin{equation*}
y(t)=[y_0(t),\dots,y_{d-1}(t)]^\top \approx [h(z_0,t),\dots,h(z_{d-1},t)]^\top.
\end{equation*}
The value $h(z_{d},t)$ is fixed for the entire time interval, it is not included in the vector function $y$.
Following \cite{Berardi_Girardi_CNSNS_2024}, we apply FD to equation \eqref{eq:state} and we obtain that its solution can be approximated by the solution of a $d-$dimensional ODE system:
\begin{align}\label{state-disc}
    \begin{aligned}
       C(y(t)) \dot{y}(t)& = \left(A_1(y(t))+A_2(y(t))\right)y(t) + Bu(t),\\
        y(0)& = y_0,
    \end{aligned}
\end{align}
where $C(y): \mathbb{R}^d \rightarrow  \mathbb{R}^{d\times d}$ is diagonal for all $y \in \mathbb{R}^d$ with the $i-$th diagonal term given by $C(y_i(t))$, whereas $A(y):=A_1(y)+A_2(y)$ is a $d \times d$ matrix for all $y \in \mathbb{R}^d$, where
$A_1(y) = \begin{pmatrix}
    \underline{0}\\
    \tilde{A}_1(y)
\end{pmatrix}$
with $\underline{0}\in\R^d$ a vector of all zeros and 
$$\tilde{A}_1(z) = \frac{1}{2\Delta z^2}{\tt trdiag}([K_{i-1}+K_i, -(K_{i+1}+2K_i+K_{i-1}), K_{i-1}+K_i],d-1),$$


\begin{equation*}
A_2(y):=\frac{1}{2 \Delta z}{\tt diag}\begin{pmatrix}
0\\ (K_0+K_1)y_0/ \Delta z +(K_2-K_0)-2 \Delta z S_1\\ (K_3-K_1)-2 \Delta z S_2\\
(K_4-K_2) -2 \Delta z S_3\\ 
\vdots \\ (K_{d}+K_{d-1})h_{N}/ \Delta z+(K_d-K_{d-2})-2\Delta z S_{d-1}.
\end{pmatrix}.
\end{equation*}
%
\first{The notation {\tt trdiag}($[v_1, v_2, v_3], d$) denotes a tridiagonal $d \times d$ matrix, where $v_1$, $v_2$, and $v_3$ specify the lower, main, and upper diagonals, respectively,}
whereas {\tt diag} maps a vector into a diagonal matrix with diagonal equal to the input vector.
Further, we denote by $K_i:=K(y_i(t))$,  $S_i:=S(y_i(t)).$  The matrix $A_2(y)$
accounts also for the root water uptake and boundary conditions in $\eqref{eq:state}$. We refer to e.g. \cite{Berardi_Girardi_CNSNS_2024,Tocci1997} for a more detailed derivation of this numerical scheme. 

\begin{remark}
\first{Due to the structure of the continuous problem \eqref{eq:state}, the first row of the matrix $A(y)$ is equal to zero. Indeed, this reflects the top boundary term which enters into the dynamical system only through the control, therefore only the matrix $B$ is affected by the dynamic boundary control. 
This peculiar structure of the matrix $A(y)$ with the first row equal to zero, leads to a singular matrix. This problem can be avoided exploiting the particular structure of the semilinear form $A(y)y$. Indeed, given a nonlinear dynamics $f(y)$, there exist infinitely many semilinear forms satisfying $f(y)=A(y)y$ in dimension greater than one. For instance, let us consider a semilinear form $A_0(x)$ satisfying $A_0(y)y = f(y)$. Then, the sum $A_0(y)+A^*(y)$ still is a semilinear form for the dynamics, provided that $A^*(y)y=\underline{0} \in \mathbb{R}^d$ for all $y \in \mathbb{R}^d$. To fill the first row, we can consider a matrix $A^*(y)$ defined as follows:
%
$$
A^*_{i,j}(x) \udef
\begin{cases}
-x_2, & \textrm{ if } i=1,\, j=1, \\
x_{j-1}-x_j,  & \textrm{if } i=1,\, j \in \{2,\ldots, d-2\}, \\
-x_{d-2}, & \textrm{if\quad}i=1,\, j=d-1, \\
0, & \textrm{otherwise}.
\end{cases}
$$
Finally, it is easy to check that $A^*(y)y=\underline{0}$ for all $y \in \mathbb{R}^d$. For more details on suitable choices for the semilinear forms of the SDRE we refer to \cite{DOLGOV2022510}.}
\end{remark}

The initial condition $y_0\in\mathbb{R}^d$ in \eqref{state-disc} is such that $(y_0)_i = h_0(z_i)$ and $B\in\mathbb{R}^d$ with $B_1 = 1$ and $B_j=0$ for $j=2,\ldots, d$ since we consider a dynamic boundary control. This setting is extremely useful to set \eqref{state-disc} and use the technique explained in the next subsection. Note that the dynamic boundary control presents a condition on the derivative of $h$ and this is reflected in \eqref{state-disc} since in the first equation, related to the boundary $z=0$, it appears also the control.

The cost functional related to the semidiscrete problem \eqref{state-disc} will be of the following form:
\begin{equation}\label{eq:cost_disc}
    \tilde{J}(u) = \int_0^T y(t)^TQ(y)y(t) + \lambda |u(t)|^2\,dt
\end{equation}
where the matrix $Q(y):\mathbb{R}^d\rightarrow \mathbb{R}^{d\times d}$ is a positive definite matrix  for all $y \in \mathbb{R}^d$, whereas $\lambda>0.$ Note that the construction of the matrix $Q(y)$ will be explained in Section \ref{sec:sdre} since it depends on the peculiar uptake function $S(h)$ due to the nonlinearity of our continuous cost \eqref{eq:cost}.

To conclude, we summarize our discretized control problem as follows: 
\begin{equation}\label{ocp_disc}
    \min_{u\in\R} \tilde{J}(u) \mbox{ such that } h \mbox{ solves } \eqref{state-disc} .
\end{equation}

\subsection{State Dependent Riccati Equation}\label{sec:sdre}
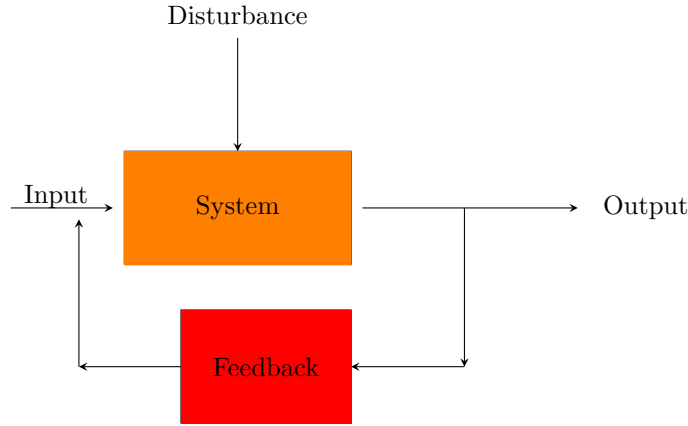
\begin{figure}[htbp]
\centering
\begin{tikzpicture}[scale=3]
    \draw (-0.3,0.3) node[centered ] {Input};
    \draw [-stealth](-0.5,0.25) -- (-0.05,0.25);
    
    \draw (0,0) -- (1,0) -- (1,0.5) -- (0,0.5) -- (0,0); 
    \fill[orange](0,0) -- (1,0) -- (1,0.5) -- (0,0.5) -- (0,0); 
    \draw (0.5,0.25) node[centered ] {System};
    
    \draw [-stealth](1.05,0.25) -- (2,0.25);
    \draw (2.3,0.25) node[centered ] {Output};
    
    \draw [-stealth](1.5,0.25) -- (1.5,-0.45);
    \draw [-stealth](1.5,-0.45) -- (1,-0.45);
    \draw (0.25,-0.7) -- (1,-0.7) -- (1,-0.2) -- (0.25,-0.2) -- (0.25,-0.7); 
    \fill[red](0.25,-0.7) -- (1,-0.7) -- (1,-0.2) -- (0.25,-0.2) -- (0.25,-0.7); 
    \draw (0.625,-0.45) node[centered ] {Feedback};
    \draw [-stealth](0.25,-0.45) -- (-0.2,-0.45);
    \draw [-stealth](-0.2,-0.45) -- (-0.2,0.2);
    \draw [-stealth](0.5,1) -- (0.5,0.5);
    \draw (0.5,1.1) node[centered ] {Disturbance};
\end{tikzpicture}
\caption{Feedback control system block diagram.}
\label{fig:feedback}
\end{figure}

As already mentioned in Section \ref{sec:intro}, we are interested in computing the control in feedback form because it offers several advantages, such as explicit dependence on the current state of the system. This closed-loop control strategy is designed to handle modeling errors, exogenous disturbances, or other variations of system parameters. This makes it particularly appealing for practical applications, where dealing with noise or small variations is common. In Figure \ref{fig:feedback}, we present a schematic realization of the feedback control.

In this subsection, we will recall one of the methods that allows us to obtain feedback control to solve \eqref{ocp_disc}: the \correct{SDRE} (see e.g. \cite{BLT07,BTB00}) approach. 

It is well-known the optimal feedback control for \eqref{ocp_disc} is linked to the value function $v(y):=\inf_u \tilde{J}(u)$, which represents the cost associated to the optimal trajectory starting from the initial condition $y \in \mathbb{R}^d$, and in this case it is given by
$$u^*(t) := -\frac{1}{\lambda}B^T\nabla v(y(t)).$$
Unfortunately, the computation of the value function requires the numerical approximation of a $d$-dimensional nonlinear PDE, the HJB equation (see e.g. \cite{BCD97}). This approach suffers from the curse of dimensionality. In a very peculiar setting, namely the Linear Quadratic Regulator (LQR) problem (\cite{dorato2000linear}), it is possible to compute the value function. This problem considers a linear dynamics with the constant matrix $A(y):=A$ in \eqref{state-disc} and a quadratic cost in \eqref{eq:cost_disc} with $Q(y):=Q$. Here, the value function is quadratic and given by $v(y) = y^T\Pi y$ where $\Pi\in\mathbb{R}^{d\times d}$ is the solution of the following algebraic Riccati equation:
$$  A^\top\Pi+\Pi A-\frac{1}{\lambda}\Pi BB^\top\Pi+Q=0.$$ This approach only relates to linear quadratic problems. We refer to e.g. \cite{Simoncini.survey.16} for efficient numerical methods for the large-scale case. Nevertheless, it is possible to extend the LQR approach to nonlinear dynamics and cost functionals provided that we deal with the dynamics given in \eqref{state-disc} and the cost \eqref{eq:cost_disc} and assuming that the value function is approximated locally by a quadratic function such that $v(y)\approx y^T\Pi(y)y$ and $\Pi(y)$ is a $d \times d$ matrix for all $y \in \mathbb{R}^d$. The unknown matrix $\Pi(y)$ is the solution of an Algebraic Riccati Equation which reads
\begin{equation}\label{sdre}
    A^\top(y)\Pi(y)+\Pi(y) A(y)-\frac{1}{\lambda}\Pi(y)B(y)B^\top(y)\Pi(y)+Q(y)=0
\end{equation}
and the feedback control, which will be suboptimal, has the following form:
\begin{equation}
u^*(t) \approx -\frac{1}{\lambda}B^T\Pi(y(t)) y(t).
\label{sdre_control}
\end{equation}
This formulation, although suboptimal since it does not solve directly the optimal control problem, corresponds to
the asymptotic stabilization of the nonlinear dynamics towards a desired configuration under suitable assumptions on $A(y),B(y),Q(y)$ and $\lambda.$
\second{To this aim, let us introduce the definitions of stabilizability and detectability.

\begin{definition}
     The pair $(A,B)$ is called stabilizable if there exists
a feedback matrix $K \in \mathbb{R}^{n \times m}$ such that $A -B K$ has all eigenvalues in the open left half complex plane.
\label{def:stab}
\end{definition}
\begin{definition}
    The pair $(A,C)$ is called detectable if the pair $(A^\top, C^\top)$ is stabilizable.
    \label{def:detec}
\end{definition}
It can be shown that, under assumptions of regularity, pointwise stabilizability, and detectability, the SDRE control achieves local asymptotic stability. This is formally stated in the following result from \cite{ccimen2008state}.
\begin{prop}
\label{prop:SDRE}
    Suppose $A(\cdot),B(\cdot)$ and $Q(\cdot)$ are $C^1(\mathbb{R}^d)$ matrix-valued function, and that the pairs $(A(x),B(x))$ and $(A(x),Q(x)^{1/2})$ are pointwise stabilizable and detectable, respectively, for each $x$, then the SDRE
method produces a closed-loop solution which is locally
asymptotically stable. 
\end{prop}
}

Note that $\eqref{sdre}$ is a matrix equation which is state dependent and has to be computed at each time iteration, \first{i.e.} every time to obtain a different state $y$. This could be very expensive especially if $d$ is very large. We refer to \cite{AKS23} for an efficient approximation of the method for large scale problems.

To summarize, the SDRE approach continuously adapts the control law to the current state of the nonlinear system, providing a powerful method for dealing with nonlinear control problems by leveraging the well-established LQR framework. The method works upon three stages: (i) linearization around the current state, (ii) solution of the algebraic Riccati equation \eqref{sdre} for the linearized problem and (iii) computation of the feedback law \eqref{sdre_control}.

\medskip

\begin{figure}[htbp]
\begin{tikzpicture}[scale=1.5]

\draw[->] (-1,0) -- (6,0) node[right] {$h$};
\draw[->] (0,-0.5) -- (0,2) node[above] {$R(h)$};

\draw (0,1) -- (-0.1,1) node[left] {$1$};
\draw (1,0) -- (1,-0.1) node[below] {$h_4$};
\draw (2,0) -- (2,-0.1) node[below] {$h_3$};
\draw (3,0) -- (3,-0.1) node[below] {$h_2$};
\draw (5,0) -- (5,-0.1) node[below] {$h_1$};

\draw[thick, red] (0,0) -- (1,0);    
\draw[thick, red] (1,0) -- (2,1);    
\draw[thick, red] (2,1) -- (3,1);    
\draw[thick, red] (3,1) -- (5,0);    
\draw[thick, red] (5,0) -- (6,0);    

\fill[red] (1,0) circle (1pt);
\fill[red] (2,1) circle (1pt);
\fill[red] (3,1) circle (1pt);
\fill[red] (5,0) circle (1pt);

\end{tikzpicture}

 \caption{The typical piece-wise linear shape of a stress function $R(h)$ as in the Feddes model \eqref{eq:Feddes}.}
    \label{fig:feddes}
    \end{figure}
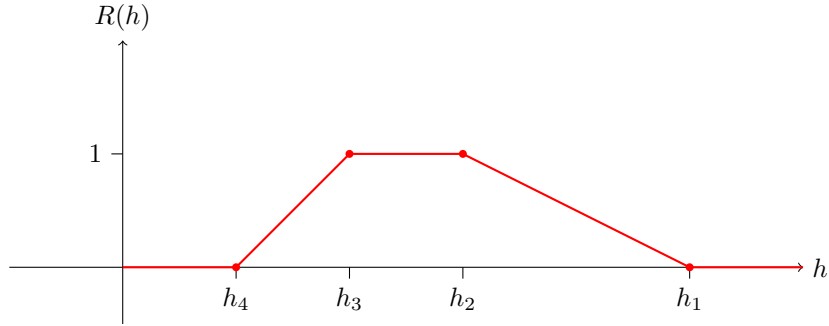

Throughout this work, we will consider the most common macroscopic root water uptake model for $S(h)$ in \eqref{eq:state}, i.e. the piece-wise linear Feddes function (\cite{Feddes_Hydrology_1976}), defined, up to a proportionality constant, as:
\begin{equation}\label{eq:Feddes}
S(h) \udef S_{max} R(h), \quad
R(h) \udef
\begin{cases}
0, & \textrm{ if }h_{1}\leq h\textrm{ or } h\leq h_{4}, \\
\frac{ h-h_{1}}{h_{2}-h_{1}}, & \textrm{if }h_{2}< h<h_{1}, \\
1, & \textrm{if\quad}h_{3}\leq  h \leq h_{2}, \\
\frac{ h-h_{4}}{h_{3}-h_{4}}, & \textrm{if }h_{4}< h<h_{3},
\end{cases}
\end{equation}
with a proper choice of parameters. A possible shape of the stress function $R(h)$ is depicted in Figure \ref{fig:feddes}. \first{Indeed, there are many smoother root water uptake functions in literature (see e.g. \cite{Li_DeJong_Coe_Ramankutty_Earth_Interactions_2006}). However, we selected the Feddes-type due to its widespread acceptance and familiarity within the soil science community. In Figure \ref{fig:Rtilde}, we also consider a regularized stress function starting from the Feddes-type.}

The function \eqref{eq:Feddes} reaches its maximum $S_{max}$ if $ h \in [h_3,h_2]$. Therefore we can rewrite the running cost in \eqref{eq:cost} as
$$
\int_{0}^T \int_0^Z |S(h(z,t)) - S_{max}|^2 dz\, dt = S_{max}^2\int_{0}^T \int_0^Z [ \chi_{(h_2,h_1)}(h)\left|\frac{ h-h_{1}}{h_{2}-h_{1}}-1\right|^2
$$
$$
+ \chi_{(h_4,h_3)}(h)\left|\frac{ h-h_{4}}{h_{3}-h_{4}}-1\right|^2 +  \chi_{(-\infty, h_4] \cup [h_1, +\infty)}(h) ] dz \, dt,
$$

where $\chi_{\omega_h}(h)$ denotes the indicator function on the set $\omega$.
Its semi-discretized version reads
$$
S_{max}^2\int_{0}^T \Delta z \sum_{i=0}^{d} [\chi_{(h_2,h_1)}(y_i(t))\left|\frac{ y_i(t)-h_{2}}{h_{2}-h_{1}}\right|^2 
$$
\begin{equation}
+ \chi_{(h_4,h_3)}(y_i(t))\left|\frac{ y_i(t)-h_{3}}{h_{3}-h_{4}
}\right|^2 +  \chi_{(-\infty, h_4] \cup [h_1, +\infty)]}(y_i(t) ) ] \, dt.
\label{eq:semidisc_cost_2}
\end{equation}
Since we are interested to recast the cost in quadratic form as in \eqref{eq:cost_disc}, it is possible to choose the state-dependent matrix $Q(y)$ as
$$
Q(y) =  S^2_{max} \Delta z  \left(Q_1(y) +Q_2(y) + Q_3(y)\right), 
$$
where
$$
Q_1(y) =  \frac{1}{|h_2-h_1|^2}(I_d+h_2^2 diag(y(t)^{-2})-2 h_2 diag(y(t)^{-1}))diag(\chi_{(h_2,h_1)}(y(t))),
$$
$$
Q_2(y) =  \frac{1}{|h_3-h_4|^2}(I_d+h_3^2 diag(y(t)^{-2})-2 h_3 diag(y(t)^{-1}))diag(\chi_{(h_4,h_3)}(y(t))),
$$
$$
Q_3(y) = diag(y(t)^{-2}) diag(\chi_{(-\infty, h_4] \cup [h_1, +\infty)]}(y(t))),
$$
where $I_d \in \mathbb{R}^{d \times d}$ is the identity matrix, $diag(y(t)^{-2})$ stands for a diagonal matrix with diagonal entry $(i,i)$ equal to $(y_i(t))^{-2}$, $diag(y(t)^{-1})$ stands for a diagonal matrix with diagonal entry $(i,i)$ equal to $(y_i(t))^{-1}$ and $diag(\chi_{\omega}(y(t)))$ stands for a diagonal matrix with diagonal entry $(i,i)$ equal to $1$ if $y_i(t) \in \omega$, otherwise it is equal to zero. It is easy to check that the expression reported in \eqref{eq:semidisc_cost_2} and the cost $\int_0^T y(t)^\top Q(y) y(t) dt$ are equivalent.

 This allows us to use the SDRE approach. We also set $\lambda = 10^{-5}$ in \eqref{eq:cost_disc}.
We note that the value $\lambda$ is very small since we do not want to penalize the control acting on the boundary to let him reach a desired configuration. \first{Furthermore, numerical evidence indicates that increasing the value of the parameter $\lambda$ does not stabilize the state within the desired interval.}

\begin{remark}\label{rmk2}
    \second{The running cost function lacks smoothness, as $R(h)$ is only Lipschitz continuous. This presents a challenge, as numerical methods for optimal control problems typically require differentiable terms. Additionally, Proposition \ref{prop:SDRE} requires that the term $Q(\cdot)$ is of class $C^1(\mathbb{R}^d)$.
    To address this limitation, 
 $R(h)$ can be approximated by a smooth function that retains its essential characteristics.}
    \second{In the case where $R(h)$ exhibits symmetry, such as $|h_4-h_3| = |h_2-h_1|$, a suitable replacement function could be of the form:}
    \second{
    \begin{equation}
    \tilde{R}(h) = e^{-(ax+b)^4},
    \label{eq:Rtilde}
        \end{equation}
    }
    \second{
    where the parameters $a$ and $b$ can be determined by solving a least-squares problem. For example, fixing $(h_1,h_2,h_3,h_4)=(0,-30,-50,-80)\,cm$, choosing $a = 0.037$ and $b=1.5$ yields the shape depicted in Figure \ref{fig:Rtilde}. In the numerical tests, in general we employ $R(h)$ as stress function and we compare the results using the two different stress functions in the first numerical test.}

    \begin{figure}[htbp]
    \centering
       \includegraphics[scale=0.5]{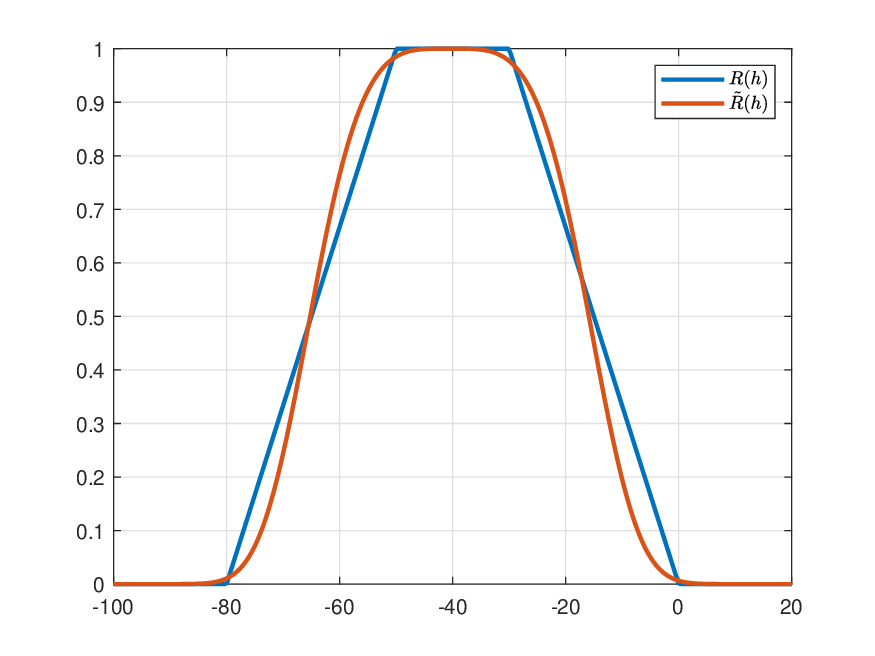}
    \caption{\second{Comparison between the stress function $R(h)$ given in \eqref{fig:feddes} and its regularized version $\tilde{R}(h)$ given in \eqref{eq:Rtilde}.}}
    \label{fig:Rtilde}
\end{figure}
\end{remark}


\section{Numerical Tests}\label{sec:test}

In this section, we present our numerical simulations. We consider two different  hydraulic models, i.e., two different choices for the water retention function $\theta(h)$ and for the unsaturated hydraulic conductivity function $K(h)$  in \eqref{eq:Richards}: the Haverkamp model and the Gardner model.
With regard to the root water uptake function $S(h)$ in \eqref{eq:Richards}, we are going to use the Feddes model \eqref{eq:Feddes}, as discussed in Section \ref{sec:sdre}, with the following choice of parameters, in cm: $h_1 = 0,\ h_2 = -30,\  h_3 = -50,\ h_4=-80$. Also, we set
$S_{max}=0.01/Z= 1.25 \cdot 10^{-4}$, where $Z=80$ cm is the soil depth.\\

We study the role of the control synthesised by the SDRE solver in these irrigation models and its stability under external perturbations. We define the mean root water uptake as
\begin{equation}
\overline{S} (t) = \frac{1}{d+1}\sum_{i=0}^{d} S(h(z_i,t)).
\label{mean_S}
\end{equation}
Whenever the mean root water uptake, at a given time $t,$ is close to the maximal value $S_{max}$, we can conclude that the water uptake is maximised for almost all the space interval at that specific time. We will consider the following the time interval \first{$[0,T]$, with $T = 1000s$}, and the space interval $[0,Z]$ is discretized with $31$ grid points. The time integration is approximated via the Matlab function \texttt{ode15s} for the resolution of stiff differential equations.
To further validate our feedback control approach we also study the problem under disturbances. This will be addressed in Section \ref{test2} and Section \ref{sec:noise2}.

The numerical simulations reported in this paper are performed on a Dell XPS 13 with Intel Core i7, 2.8GHz and 16GB RAM. The codes are written in Matlab R2023b.


\subsection{Test 1: Haverkamp model}
In the first test, we consider the Haverkamp model where the water retention curve and hydraulic conductivity are defined as follows:
\begin{equation}\label{eq:Haverkamp}
\theta(h) = \frac{\alpha \left(\theta_S - \theta_r \right)}{\alpha + \abs{h}^{\beta_2}} + \theta_r; \quad
K(h) = K_S  \frac{A}{A+\abs{h}^{\beta_1}}.
\end{equation}
We remark that some parameters have a physical meaning: $\theta_S$ represents the saturated water content, $\theta_r$ is the residual water content, with the obvious constraint that $\theta_r \leq \theta \leq \theta_S$; $K_S$ represents the saturated hydraulic conductivity. Other parameters in the following formulas are fitting parameters from real data, in practice. 
In this framework, we consider the following set of parameters, \first{as in \cite{Celia_et_al}}: 
$$K_S = 34 \mbox{cm/h},\; A = 1.175\cdot10^{6}, \; \alpha = 1.611\cdot10^{6},$$ 
$$\theta_S = 0.287, \; \theta_r = 0.075,\; ,\;\beta_1 = 4.74, \; \beta_2 = 3.96,$$
with the same initial and boundary conditions: we fix a constant initial condition $h_0(z) = -61.5$ cm and we use the same value for the bottom boundary condition $h_B(t) = -61.5$cm, 
while we consider as initial condition for the top boundary condition $h_T = -20.73$ cm. 

\medskip 

\begin{figure}[htbp]
    \centering
       \includegraphics[scale=0.4]{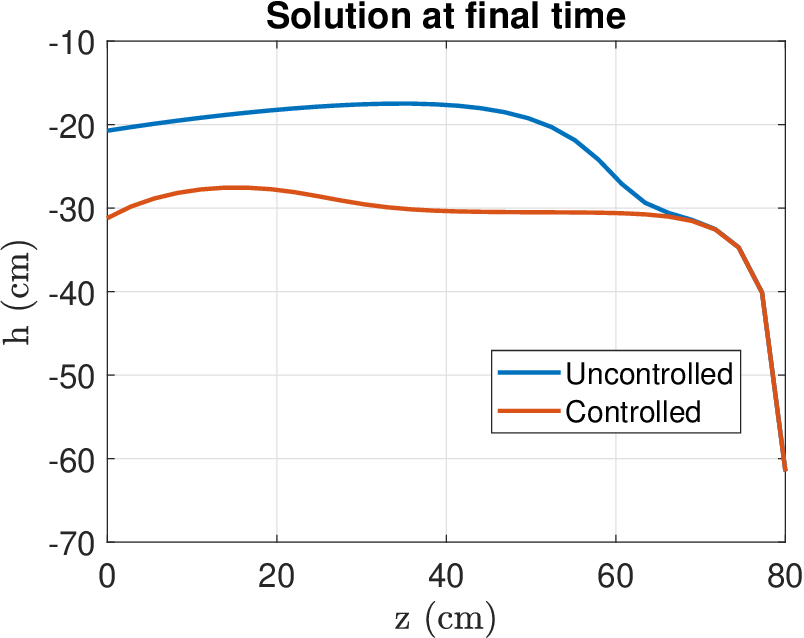}
     \includegraphics[scale=0.4]{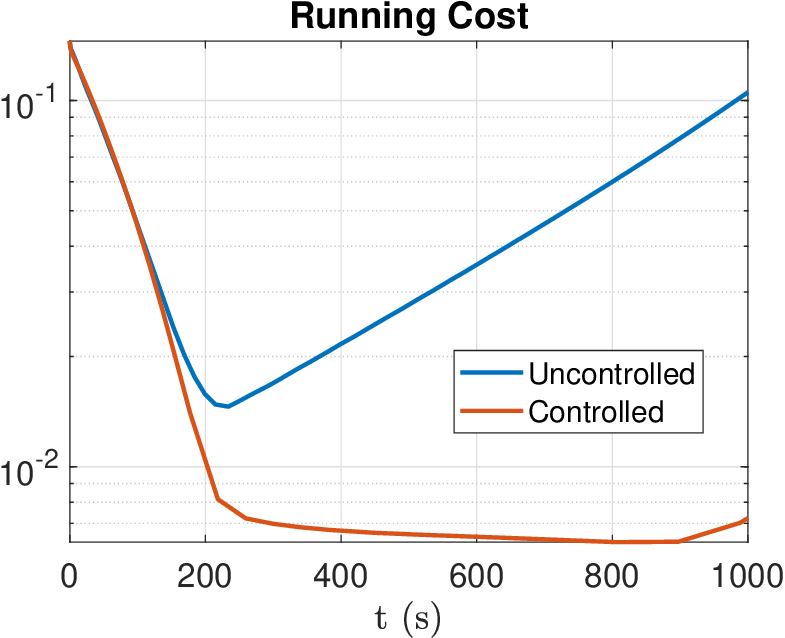}\\
      \includegraphics[scale=0.4]{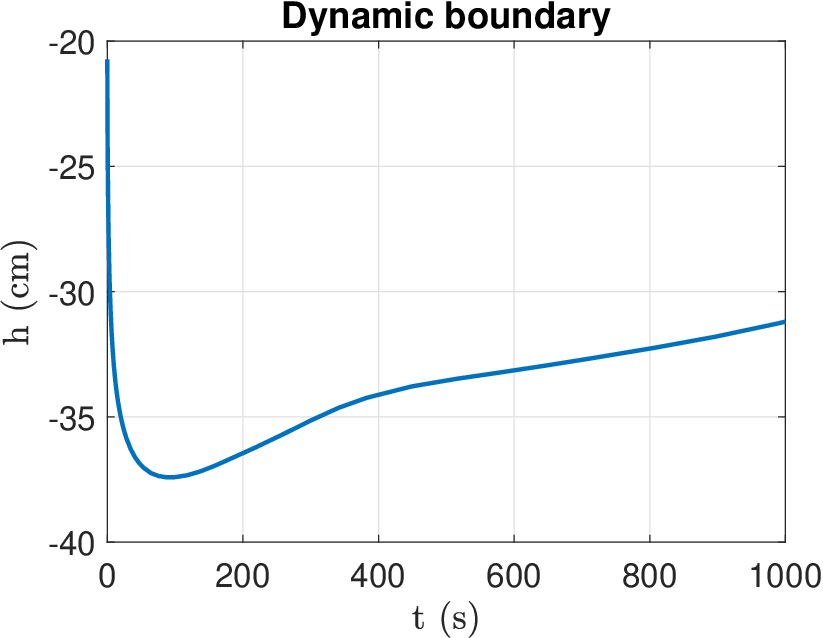}
            \includegraphics[scale=0.4]{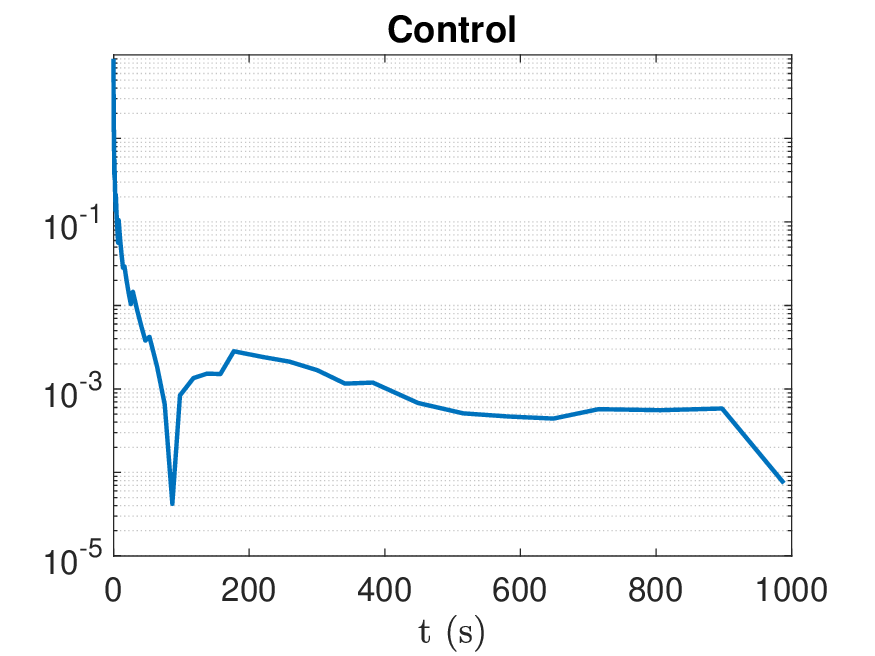}
    \caption{Test 1. Solution at final time \first{$T=1000s$} for the controlled and uncontrolled solutions (top-left), running cost (top-right), the dynamic boundary condition (bottom-left) and control (bottom-right).}
    \label{fig1:sol}
\end{figure}

In the top-left panel of Figure \ref{fig1:sol}, we compare the solution at final time $T=1000s$ for the controlled and uncontrolled solution. One may see that even if the control is acting only as a left boundary condition, the controlled solution at final time stays closer to the desired interval $[-60,-30]$ with respect to the uncontrolled one.
By the comparison of the running cost displayed in the top-right panel of Figure \ref{fig1:sol} one can verify that this closeness to the desired interval is global in time, since the cost functional related to the controlled solutions is always \first{less than or equal to the cost in the uncontrolled case}. In the bottom panel of Figure \ref{fig1:sol}, we report the computed control (right panel) and the corresponding left dynamic boundary condition (left panel). We immediately note that the control starts with high values for the first time instances and its magnitude is subject to a fast decrease. This is reflected in the behaviour of the dynamic boundary, passing from $-20$ cm to $\approx -35$ cm in the very first time instances, while for the rest of the interval it starts increasing with a small time derivative. Indeed, the control is very active at the beginning since in the left boundary we do no start close to the desired configuration.

\begin{figure}[htbp]
    \centering
    \includegraphics[scale=0.6]{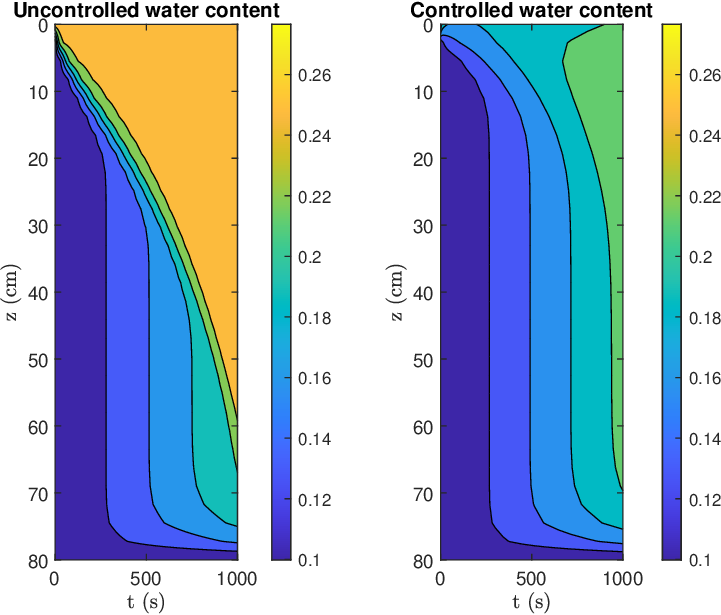}
    \caption{Test 1. Representation of $\theta(h)$ given in \eqref{eq:Haverkamp} for the uncontrolled solution (left) and controlled solution (right). The $x-$axis represents the temporal coordinates whereas the $y-$axis the spatial domain.}
    \label{fig1:theta}
\end{figure}

\begin{figure}[htbp]
    \centering
    \includegraphics[scale=0.45]{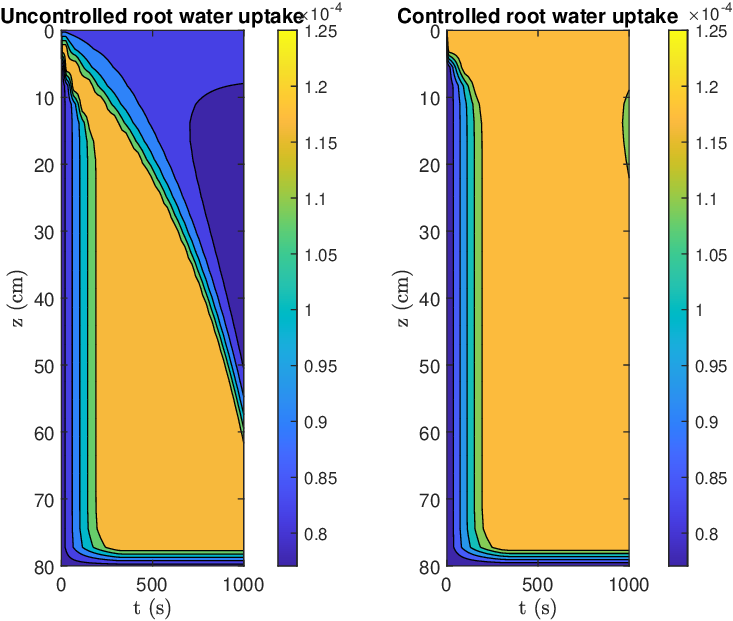}
    \includegraphics[scale=0.4]{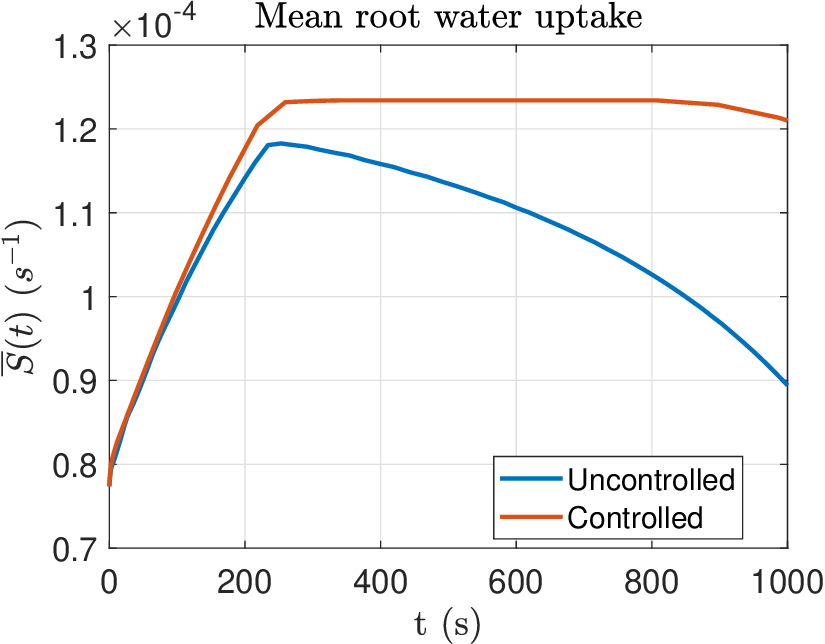}
    \caption{Test 1. Left: Contour lines of $S(h)$ given in \eqref{eq:Feddes} for the uncontrolled and controlled solutions. Right: mean root uptake $\overline{S}(t)$ for the uncontrolled and controlled solutions. }
    \label{fig1:S}
\end{figure}

In Figure \ref{fig1:theta}, we compare the water content $\theta(h)$ given in \eqref{eq:Haverkamp} for the uncontrolled solution (left panel) and the controlled solution (right panel). We note that the control does its job, allowing to save water (as can be seen from the comparison of the controlled and uncontrolled mode in the same plot).
Finally, in Figure \ref{fig1:S} we compare the root water uptake function for the uncontrolled and the controlled case through the contours of the function $S(t,z)$ (left and middle panel). It is also interesting to note how the action of the bottom boundary condition influences the solution in the first time instances and, afterwards, only in a small region close to that boundary. Clearly the control cannot act there. For the sake of completeness, we also report their means $\overline{S}(t)$ introduced in \eqref{mean_S} in the right panel of Figure \ref{fig1:S}. We recall that the maximum value $S_{max}$ for the uptake function is $1.25 \cdot 10^{-4}$ and we note that the controlled mean root water uptake stays close to this value after $\approx 250s$, while in the same time interval the uncontrolled one decreases.  

    

\begin{figure}[htbp]
    \centering
       \includegraphics[scale=0.4]{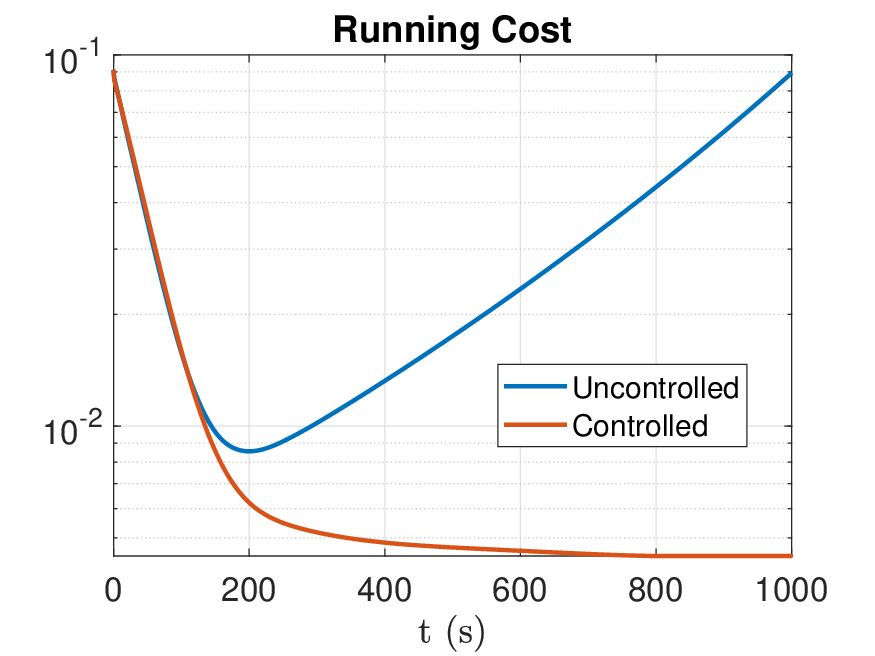}
     \includegraphics[scale=0.4]{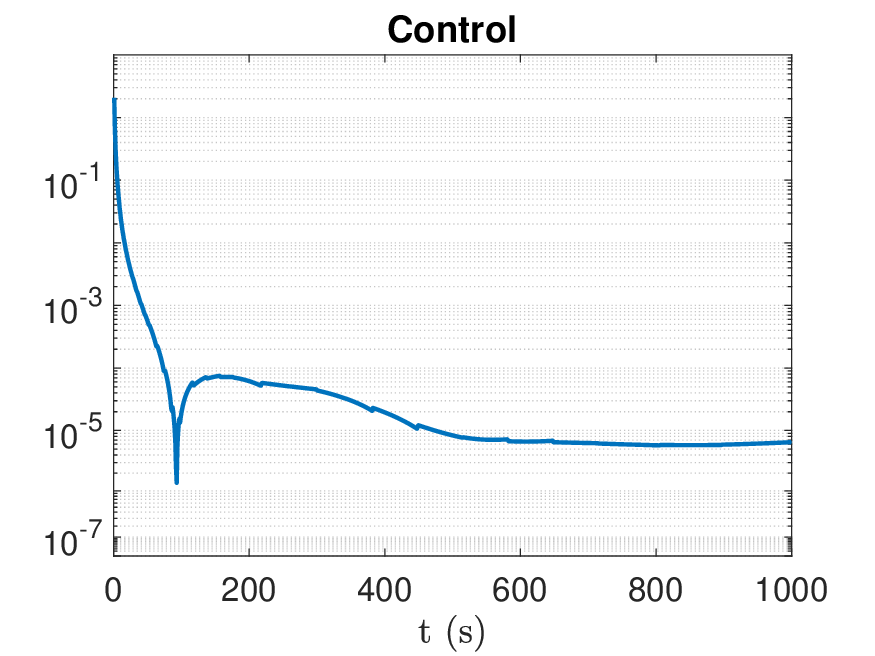}
    \caption{\second{Test 1. Running cost (left) and control (right) with stress function $\tilde{R}(h)$.}}
    \label{fig:sol_Rtilde}
\end{figure}

\begin{figure}[htbp]
    \centering
       \includegraphics[scale=0.4]{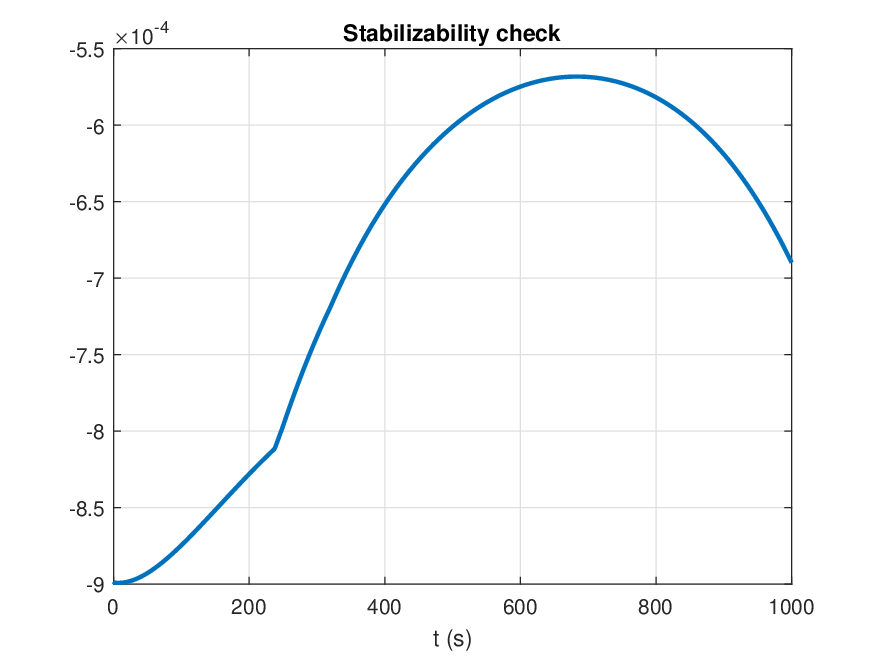}
     \includegraphics[scale=0.4]{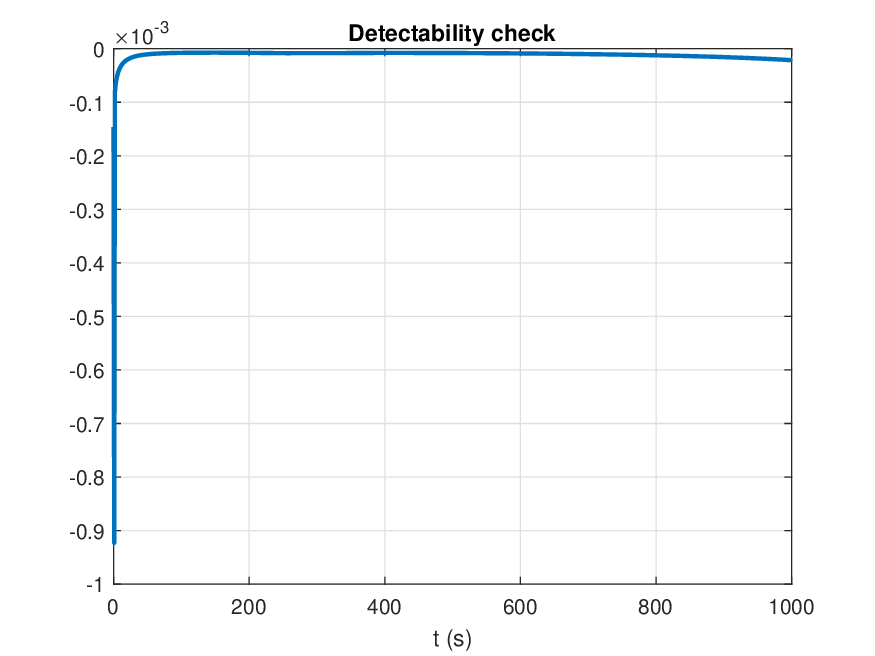}
    \caption{\second{Test 1. Maximum real part of the eigenvalues of $A(x)-BK(x)$ (left) and $A(x)^\top - 100(Q(x)^{1/2})^\top$ (right) over the time integration interval.}}
    \label{fig:stab_check}
\end{figure}

\second{In Remark \ref{rmk2}, we observed that the control law lacks smoothness, and the running cost of the controlled dynamics begins to increase toward the end of the time interval. This behavior may be attributable to the lack of smoothness in the running cost. To address this, we replaced the stress function $R(h)$ with its regularized version shown in \eqref{eq:Rtilde}. The resulting behavior is displayed in Figure \ref{fig:sol_Rtilde}. In this case, the running cost demonstrates a non-increasing trend, and the control signal appears smoother, suggesting that the enhanced smoothness of the running cost contributes to achieving a more refined and stable control outcome. Then, we have verified if the hypothesis of stabilizability and detectability introduced in Definition \ref{def:stab} and \ref{def:detec} are satisfied. In Figure \ref{fig:stab_check}, we display the maximum real part of the eigenvalues of the closed loop matrix $A(x)-BK(x)$ in the left panel and $A(x)^\top - 100(Q(x)^{1/2})^\top$ in the right panel. We note that they remain negative for the entire time interval, indicating that the hypothesis of Proposition \ref{prop:SDRE} are satisfied and hence the SDRE solution is locally asymptotically stable.}

\subsection{Test 2: Haverkamp model with noise}
\label{test2}
In this test case, we consider the same setting of the previous test and we investigate the role of the feedback control in the stability under external perturbations.

\begin{figure}[htbp]
    \centering
       \includegraphics[scale=0.4]{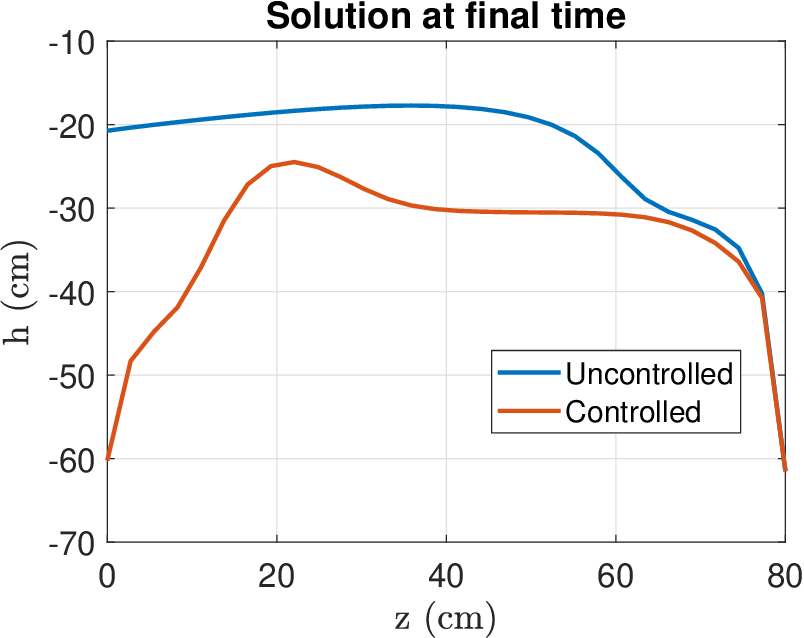}
     \includegraphics[scale=0.4]{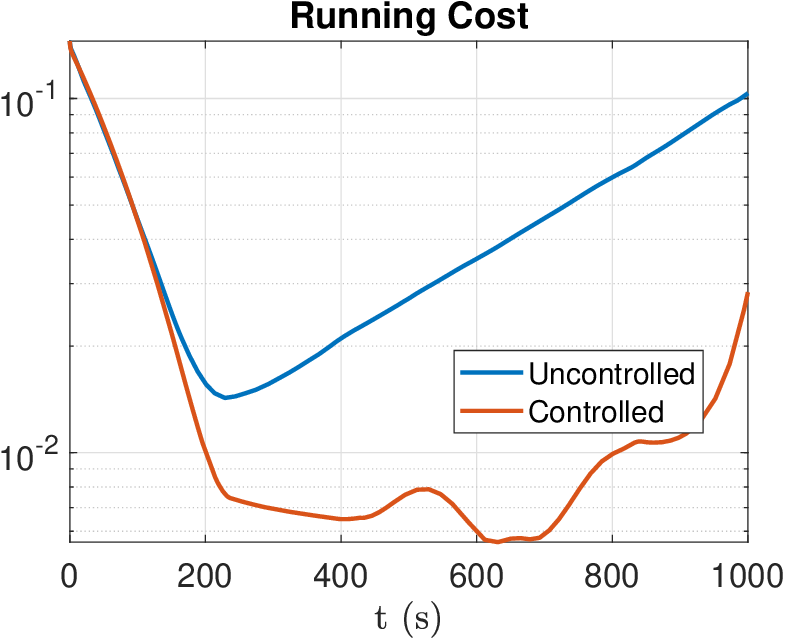}\\
      \includegraphics[scale=0.4]{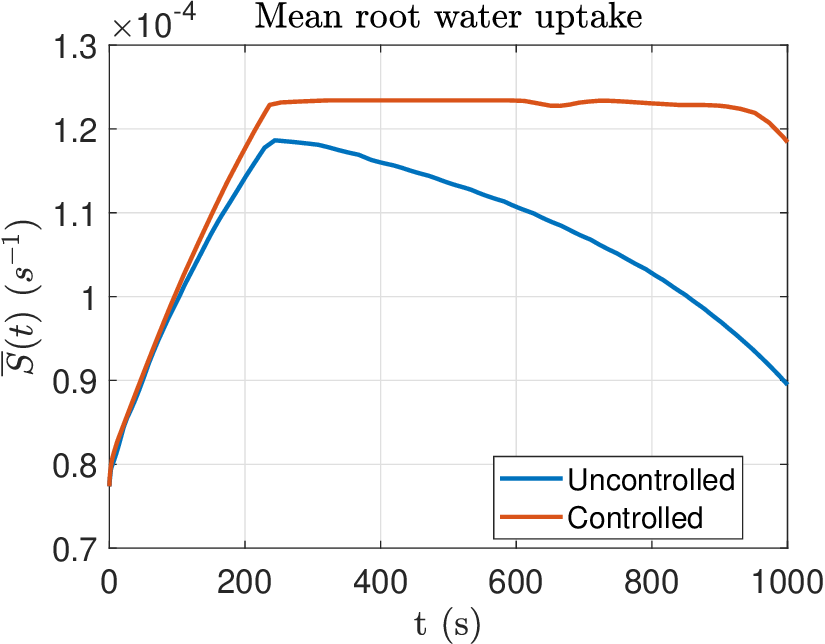}
      \includegraphics[scale=0.4]{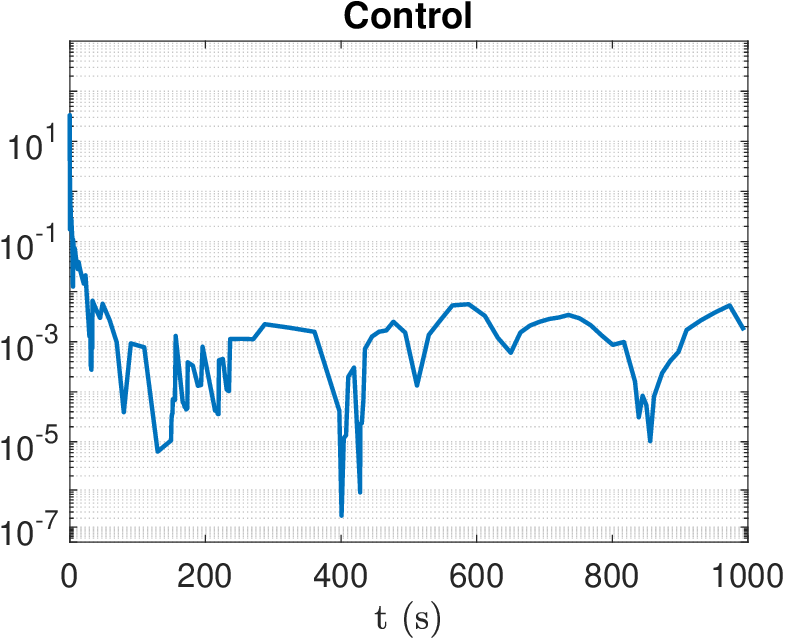}
    \caption{Test 2. Solution at final time \first{$T=1000s$} for the controlled and uncontrolled solutions (top-left), running cost (top-right), mean root water uptake (bottom-left) and the control (bottom-right) in presence of a noisy term.}
    \label{fig3:sol}
\end{figure}

 In particular we suppose that the term $K(h)$ is affected by a source of noise or uncertainty, substituting it with $K_{noisy}(h) := K(h)(1+\epsilon \eta)$, with $\eta$ a uniform random variable taking value in the interval $[0,1]$ and $\epsilon$ represent the amplitude of the noise. Let us fix $\epsilon = 10^{-5}$. The optimal control problem \eqref{ocp} is now solved with hydraulic conductivity $K_{noisy}(h)$. We point out that $K_{noisy}(h)$ is a random variable and hence its value will vary along the time integration of the \first{Richards'} equation. In Figure \ref{fig3:sol}, we display the results obtained for this test case. The uncontrolled solution does not present any visual difference with respect to the case without noise. On the other hand, we note that the control strategy is highly influenced by the noise, especially observing the control behaviour in the bottom-right panel of Figure \ref{fig3:sol}. Indeed, the noise introduces many oscillations on the control, resulting in a irregular behaviour for the controlled solution at final time and running cost (top-left and top-right panels). \first{The presence of the noise has slightly changed the results of the running cost compared to Figure \ref{fig1:sol}. Indeed here, for the first time instances the controlled case has a running cost a bit larger than the uncontrolled case. However, the whole cost for the controlled case is smaller than the uncontrolled.} Nevertheless, the controlled strategy is still able to reduce the cost functional and to keep the mean root water uptake close to the value $S_{max}$, as shown in the bottom-left panel.

Finally in Table \ref{table1}, we compare the total costs obtained integrating the uncontrolled and the controlled running costs along the time interval \first{$[0,T]$}. In both cases the ratio between the uncontrolled and the controlled costs is almost 3, while we see a slight increase of the cost for the controlled trajectory in Test 2 due to the presence of the noise and the nonlinear nature of the problem. As already mentioned, we did not observe substantial differences in the uncontrolled solution. Indeed, the total cost is very similar (first row Table \ref{table1}).
\begin{table}[hbht]
\centering
\begin{tabular}{c|c|c}     
 & Test 1 & Test 2  \\ \hline

Uncontrolled & 45.68  & 45.23  \\ 
Controlled &  16.37  & 17.23 \\  \hline

 \end{tabular}
  \caption{Total cost for the uncontrolled and controlled trajectories for the different tests.}
 \label{table1}
\end{table}

\subsection{Test 3: Gardner model}

In the third test, we consider the Gardner model where the water retention function and conductivity function are defined respectively as

\begin{equation}
\label{eq:Gardner}
\theta(h) \udef \theta_{r}+(\theta_{S}-\theta_{r})e^{\rho h}, \qquad K\left( h \right)\udef K_S e^{\rho h}.
\end{equation}

We refer to e.g.  \cite{Suk_Park_J_Hydrology_2019,Berardi_Difonzo_JCD_2022} for more details on the Gardner hydraulic functions. The definition of these functions usually depends on the sign of the variable $h$, but since we are considering problems for which $h$ is strictly negative, we report just the definitions for \first{this case}. We set $\rho = 0.1 \,cm^{-1}$,  $K_S = 1 \, cm/s, \theta_S = 0.48$ and $\theta_{r} = 0$, \first{as reported in \cite{Berardi_Difonzo_JCD_2022}}. 

Given the functions $\theta(h)$, $C(h)$ can be obtained in the following way
$$
C(h) = \frac{\partial \theta }{\partial h} = \rho  (\theta_{S}-\theta_{r})e^{\rho h}.
$$

\begin{figure}[htbp]
    \centering
       \includegraphics[scale=0.4]{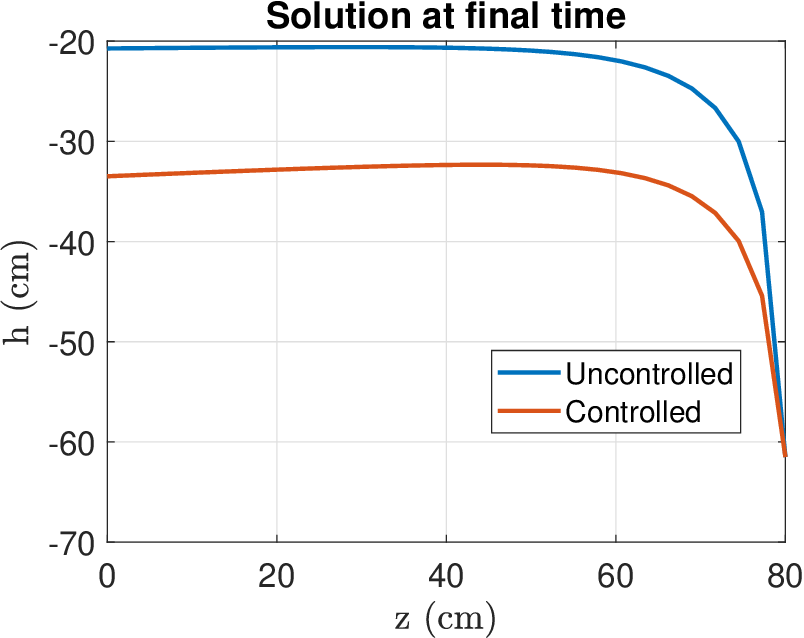}
     \includegraphics[scale=0.4]{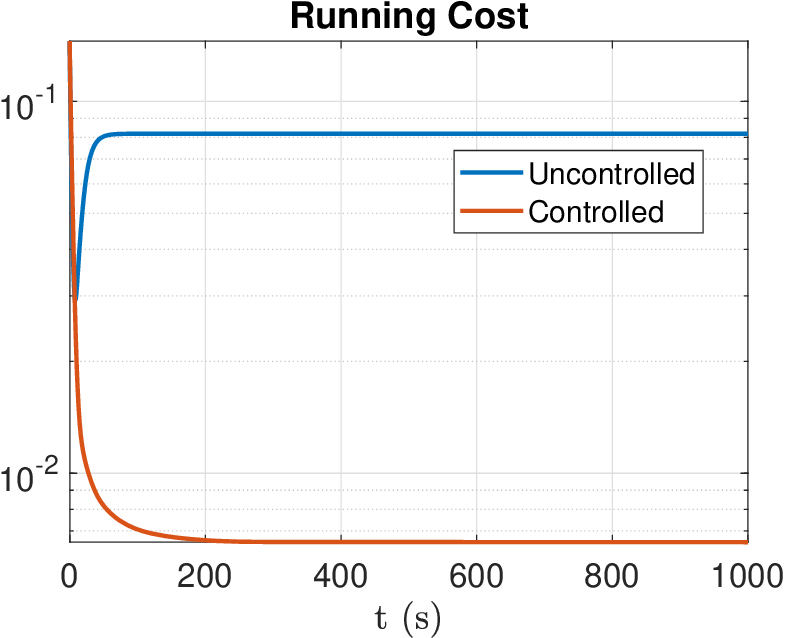}\\
      \includegraphics[scale=0.4]{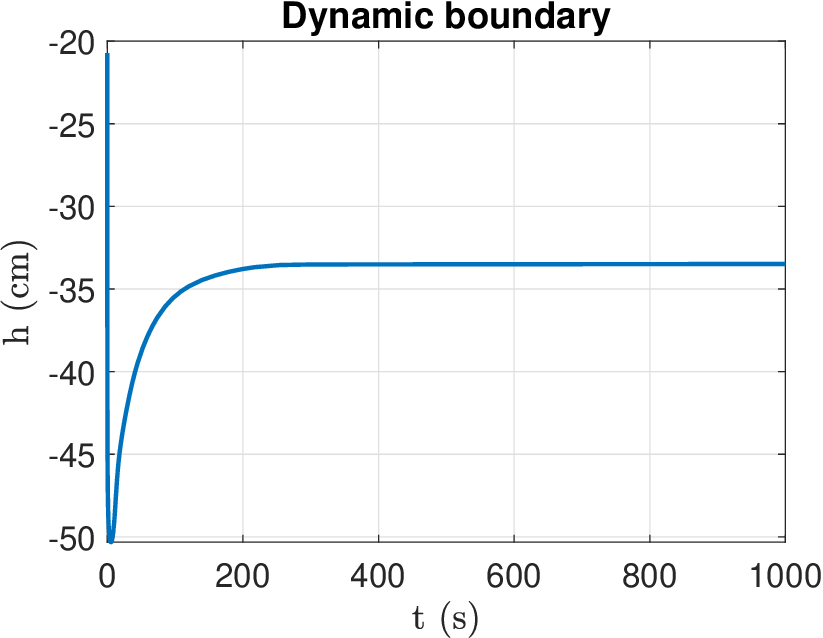}
       \includegraphics[scale=0.4]{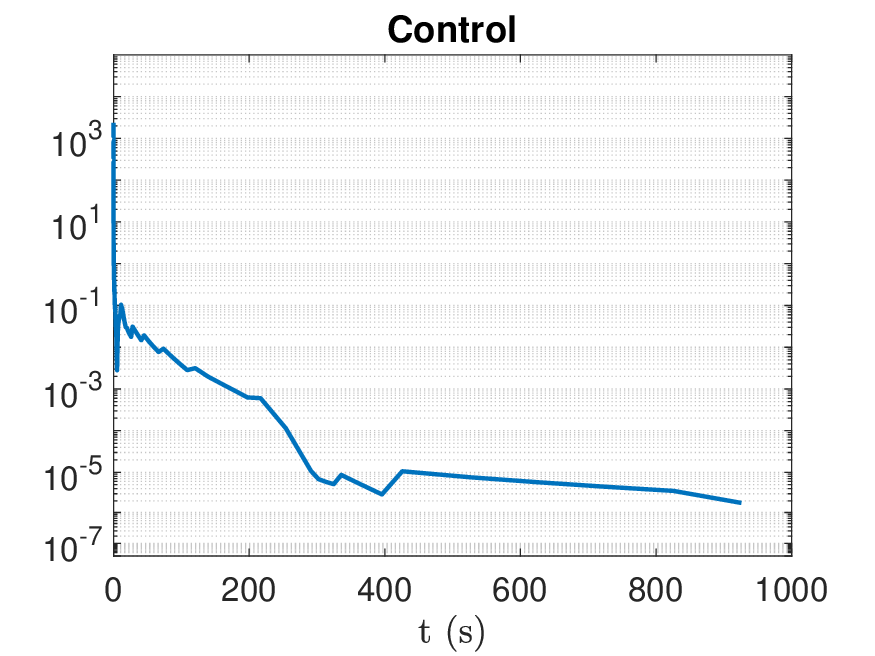}
    \caption{Test 3. Solution at final time $T=1000 \,s$ for the controlled and uncontrolled solutions (top-left), running cost (top-right), dynamic boundary condition (bottom-left) and the control (bottom-right).}
    \label{fig4:sol}
\end{figure}

In the top-left panel of Figure \ref{fig4:sol}, we display the configuration of the solutions at the final time $T=1000$s. We immediately note that almost the entire controlled solution is contained in the desired interval $[-60,-30]$, except for the bottom boundary condition which is fixed outside of the interval by the definition of the problem. In the top-right panel we show the comparison between the running costs. The cost of the controlled solution is higher just for the first time instances, due to the initial high magnitude of the control ($\approx 10^{3}$), while it stays below the uncontrolled one with one order of magnitude of difference. The initial high value of the control is reflected in the behavior of the dynamic boundary condition, which immediately goes from $-20$ cm to $-50$ cm. Afterwards, it increases slowly until it reaches a plateau with value $\approx -34$ cm, which again is included in the interval where the uptake function reaches its maximum.


\begin{figure}
    \centering
    \includegraphics[scale=0.6]{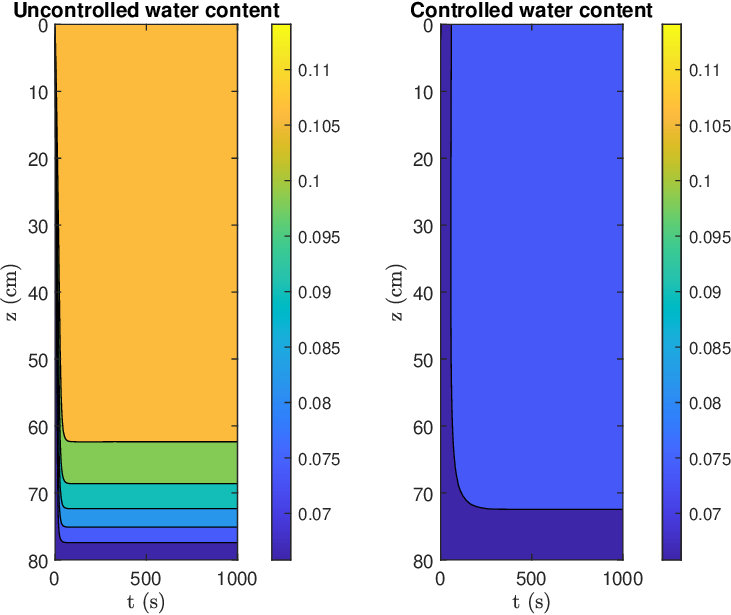}
    \caption{Test 3. Representation of $\theta(h)$ given in \eqref{eq:Gardner} for the uncontrolled solution (left) and controlled solution (right). The $x-$axis represents the temporal coordinates whereas the $y-$axis the spatial domain.}
    \label{fig4:theta}
\end{figure}

In Figure \ref{fig4:theta}, the water content $\theta(h)$ given in \eqref{eq:Gardner} for the uncontrolled solution (left panel) and the controlled solution (right panel) are compared. One can see that the control procedure allows to save a significant amount of water, yet ensuring the desired root uptake.

\begin{figure}[htbp]
    \centering
    \includegraphics[scale=0.45]{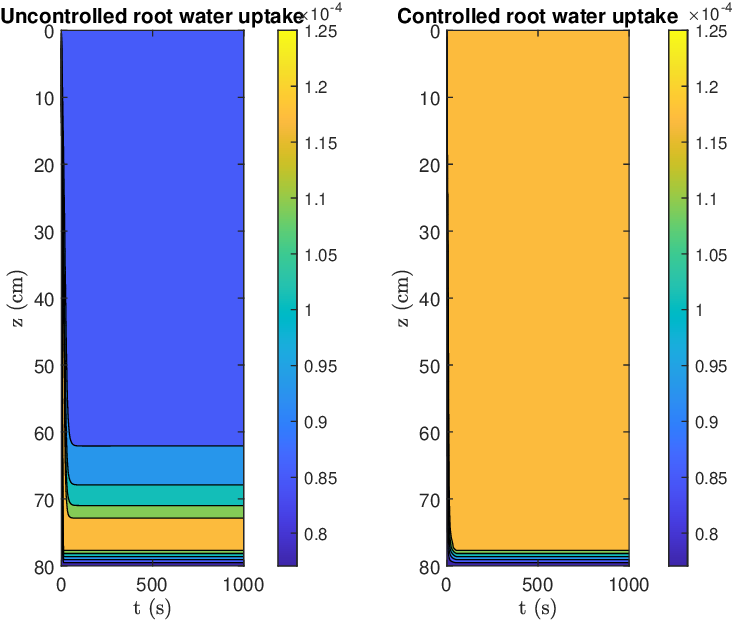}
    \includegraphics[scale=0.4]{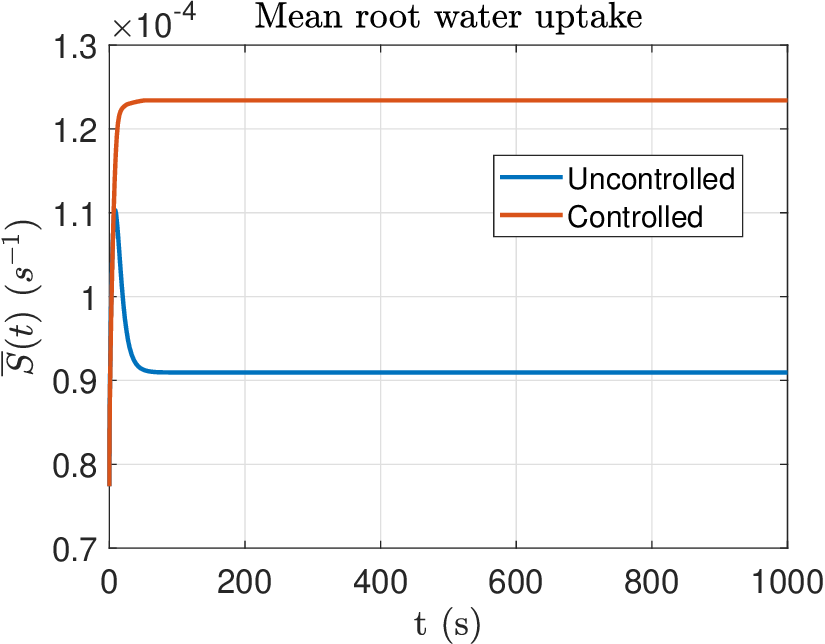}
    \caption{Test 3. Left: Representation of $S(h)$ given in \eqref{eq:Feddes} for the uncontrolled and controlled solutions. Right: mean root uptake $\overline{S}(t)$ for the uncontrolled and controlled solutions.}
    \label{fig4:S}
\end{figure}

Finally, in Figure \ref{fig4:S} we show the uncontrolled and the controlled root water uptake function (left and middle panel) and their means $\overline{S}(t)$ in the right panel. On one hand the mean root water uptake reaches almost the maximum value $S_{max}$ in the first time instances, on the other hand also the uncontrolled one reaches immediately an equilibrium, but to a lower value ($\approx 9 \cdot 10^{-5}$).

\subsection{Test 4: Gardner model with noise}\label{sec:noise2}

\begin{figure}[htbp]
    \centering
       \includegraphics[scale=0.4]{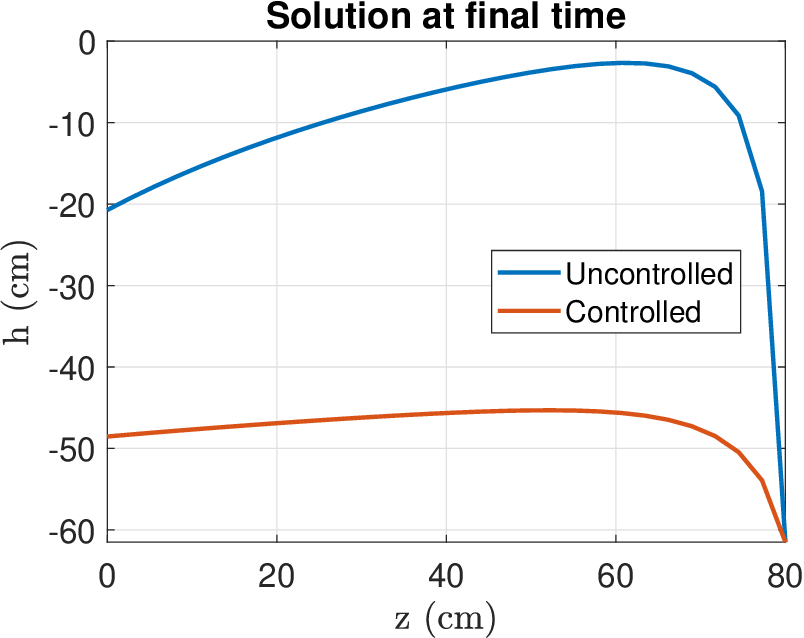}
     \includegraphics[scale=0.4]{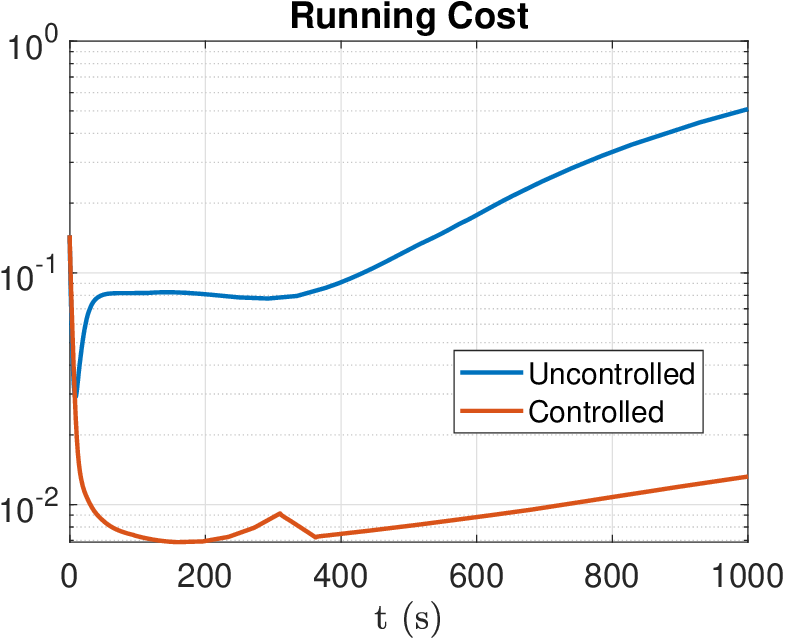}\\
      \includegraphics[scale=0.4]{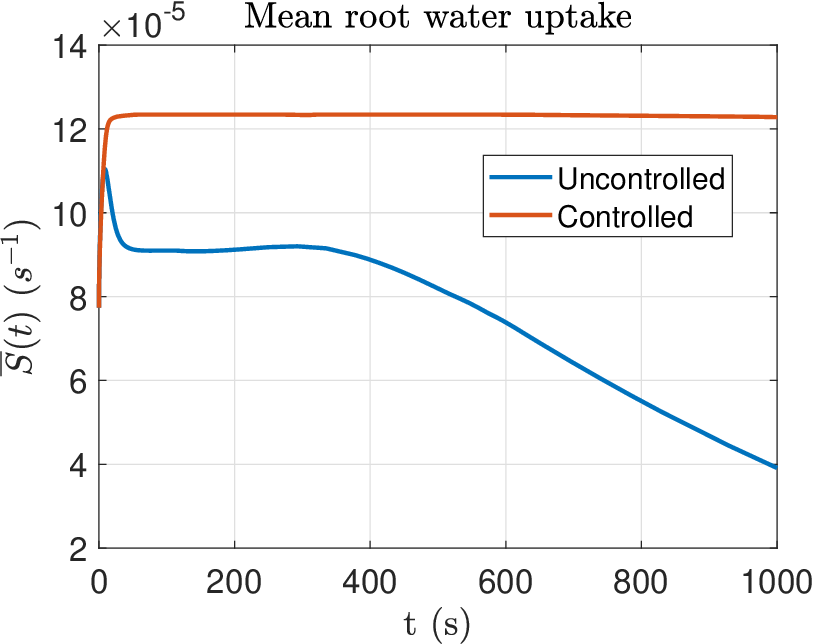}
      \includegraphics[scale=0.4]{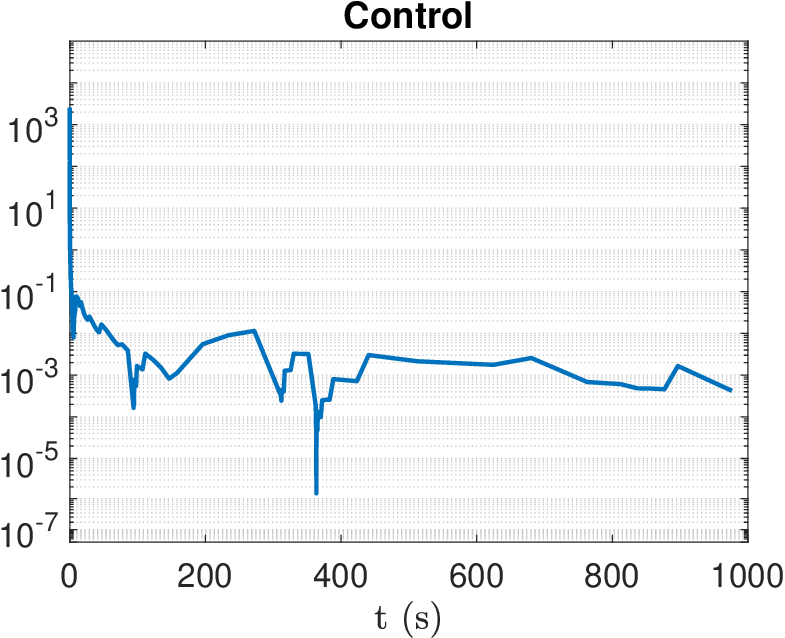}
    \caption{Test 4. Solution at final time \first{$T=1000s$} for the controlled and uncontrolled solutions (top-left), running cost (top-right), mean root water uptake (bottom-left) and control (bottom-right) in presence of a noisy term.}
    \label{fig5:sol}
\end{figure}

Last example deals with the stability of the control strategy. Again, we study the system under the influence of an external source of noise for the Gardner model. Once more, the hydraulic conductivity is perturbed considering $K_{noisy}(h) = K(h)(1+\epsilon \eta)$, where $\epsilon$ and $\eta$ have been introduced in Test \ref{test2}. We fix $\epsilon = 10^{-6}$ \first{since in this case the hydraulic conductivity is an exponential function and its perturbation affects strongly the problem}. In the top-left panel of Figure \ref{fig5:sol}, it is possible to note how far is the uncontrolled solution at final time from the configuration of the previous test in absence of noise. This results in a higher running cost and a lower mean root water uptake, since the solution is more distant from the desired interval. 
On the other hand, the controlled solution at final time keeps staying in the opportune interval, resulting in a behaviour similar to the previous case. The presence of the noise reflects again in a presence of oscillations in the control signal, which appears less smooth than the case without noise.

Finally, in Table \ref{table2} the total costs for the uncontrolled and controlled solutions for the Gardner model are reported for the last two numerical tests. In absence of noise, the SDRE controller manages to reduce the total cost with a ratio of almost $11$. In presence of noise the uncontrolled total cost increases more than the double, while the controlled is subject to a small increase, resulting in a ratio of almost $20$, showing the real strength of the feedback control. 
\begin{table}[hbht]
\center
\begin{tabular}{c|c|c}     
 & Test 3 & Test 4  \\ \hline

Uncontrolled & 80.84  &   191.71  \\ 
Controlled &  7.25  &  9.60 \\  \hline

 \end{tabular}
  \caption{Total cost for the uncontrolled and controlled trajectories for the different tests.}
 \label{table2}
\end{table}
\section{Conclusions}
In this work we have introduced a first approach to compute a feedback control in the context of Richards' equation, modeling water flow in unsaturated porous media. Specifically, we have introduced a dynamic boundary control problem which was never tackled before, aimed at maximizing water uptake by plant roots and simultaneously minimizing irrigation.
 In general, it is possible to observe that for the problems discussed, the control introduces first a strong modification of the boundary condition, revealing the suitable initial framework for the control problem. Afterwards, it is characterized by a slow and smooth change, sometimes reaching an equilibrium depending on the model considered.
Future directions will include a theoretical study of existence and uniqueness of the control problem presented and more sophisticated feedback control strategies aiming at the optimality. \first{Additionally, this approach will be applied to a lysimeter-scale test, involving plant growth in a soil-filled plot equipped with soil water content sensors}, in order to validate this approach in a real life scenario. Furthermore, it will be interesting to investigate this model including a spatial dependent root water uptake  function.

\section*{Declarations}

\paragraph{\bf Funding} A. Alla, M. Berardi, L. Saluzzi are members of the INdAM-GNCS activity
group. A. Alla is part of INdAM - GNCS Project \first{“}Metodi numerici
innovativi per equazioni di Hamilton-Jacobi” (CUP E53C23001670001). L. Saluzzi is part
of INdAM - GNCS Project \first{“}Metodi di riduzione di modello ed approssimazioni di rango basso per problemi alto-dimensionali" (CUP E53C23001670001) and \first{“}titolare di borsa per l’estero dell’Istituto Nazionale di Alta Matematica”.

The work of A.A. has been carried out within the “Data-driven discovery and control of multi-scale interacting artificial agent systems”, and received funding from the European Union Next-GenerationEU - National Recovery and Resilience Plan (NRRP) – MISSION 4 COMPONENT 2, INVESTIMENT 1.1 Fondo per il Programma Nazionale di Ricerca e Progetti di Rilevante Interesse Nazionale (PRIN) – Project 2022 No.~P2022JC95T. This manuscript reflects only the authors’ views and opinions, neither the European Union nor the European Commission can be  considered responsible for them. A.A. is also supported by MIUR with PRIN project 2022 funds (2022238YY5, entitled
"Optimal control problems: analysis, approximation ")

The work of M.B. has been carried out within the “SAFER MESH: Sustainable mAnagement oF watEr Resources: ModEls and numerical MetHods”, and received funding
from the European Union Next-GenerationEU - National Recovery and Resilience
Plan (NRRP) – MISSION 4 COMPONENT 2, INVESTIMENT 1.1 Fondo per il
Programma Nazionale di Ricerca e Progetti di Rilevante Interesse Nazionale (PRIN)
– CUP B53D23027870001. This manuscript reflects only the authors’ views and
opinions, neither the European Union nor the European Commission can be considered responsible for them. Moreover, M.B. thanks Mrs Domenica Livorti for supporting the project's activities. 


\paragraph{\bf Declaration of competing interest}
The authors declare that they have no known competing financial interests or personal relationships that could have appeared to influence the work reported in this paper.

\bibliographystyle{plain}
\bibliography{bibliogroundwater_new,bibsdre,references}

\end{document}